\documentclass[leqno]{article}
\usepackage{amssymb,amsmath, amsthm, amsfonts}
\usepackage{graphicx, epsf,epsfig, float}
\usepackage{latexsym}

\newcommand{\disp}{\displaystyle}
\newcommand{\bi}{\begin{itemize}}
\newcommand{\ei}{\end{itemize}}
\newcommand{\ba}{\begin{array}}
\newcommand{\ea}{\end{array}}
\newcommand{\benu}{\begin{enumerate}}
\newcommand{\eenu}{\end{enumerate}}
\newcommand{\be}{\begin{equation}}
\newcommand{\ee}{\end{equation}}

\graphicspath{{./jpgfiles/}}

\begin{document}
\title{Axially symmetric volume constrained anisotropic mean curvature flow}
\author{ By  B{\footnotesize BENNETT} P{\footnotesize ALMER}
 and  W{\footnotesize ENXIANG} Z{\footnotesize HU}}

\date{}
\maketitle
\newtheorem{theorem}{Theorem}[section]
\newtheorem{cor}{Corollary}[section]
\newtheorem{prop}{Proposition}[section]
\newtheorem{lemma}{Lemma}[section]
\newtheorem{condition}{Condition}[section]
\newtheorem{example}{Example}[section]

\newtheorem{definition}{Definition}[section]
\newtheorem{remark}{Remark}[section]
\newtheorem{conjecture}{Conjecture}[section]
\newtheorem{claim}{Claim}[section]
\newtheorem{question}{Question}[section]
\newcommand{\rf}[1]{\mbox{(\ref{#1})}}

\renewcommand{\thefootnote}{\fnsymbol{footnote}}
\abstract{We study the long time existence theory for a non local flow associated to a free boundary problem for a trapped non liquid drop. The drop has free boundary components on two horizontal plates and its free energy is anisotropic and axially symmetric. For axially symmetric initial surfaces with sufficiently large volume, we show that the flow exists for all time.

Numerical simulations of the curvature flow are presented.}\\[5mm]
\section{Introduction}
The evolution of interfaces of structured materials is of interest in a wide range of disiplines related to materials science \cite{JT}. Structured materials such as crystals, polycrystals and liquid crystals have a surface energy which is anisotropic; their energy density depends on the direction of the surface at each point. Over the years, various methods have been developed to track these interfaces, including the phase-field and level set methods. Here, we will consider a particular free boundary problem utilizing the anisotropic mean curvature flow. 

The mean (isotropic)  curvature flow with constrained volume was considered in \cite{H}. In relation to the free boundary problem considered here, the papers of Athanassenas \cite{Ath}, \cite{Ath2} are particularly relevant .  Also, in the recent paper \cite{crm},  volume preserving mean curvature flow in a Riemannian setting is studied.

Volume constrained anisotropic mean curvature flow for hypersurfaces was considered in
\cite{An}, \cite{Mcc}. In these papers, the emphasis is on the evolution of closed convex hypersurfaces.

Consider an anisotropic surface energy which assigns to a sufficiently smooth surface with unit normal $\nu$ the value
\begin{equation}
\label{F}{\cal F}(\Sigma)=\int_\Sigma \gamma(\nu) \:d\Sigma \:.\end{equation}
The function $\gamma:S^2\rightarrow {\bf R}^+$ is assumed to satisfy a {\it convexity condition}: the surface
\begin{equation}
\label{W} W=\partial  \bigcap_{n\in S^2} \{ Y\cdot n\le \gamma(n)\}\:,\end{equation}
which is known as the {\it Wulff shape} .
In this paper, it will be assumed that
\begin{enumerate}
\item[(W1)] $W$ is a smooth,uniformly convex surface of revolution with vertical rotation axis.
\item[(W2)] $W$ is symmetric with respect to reflection through the horizontal plane $z=0$.
\end{enumerate}
Although we will not assume it in general, the following condition will also enter into discussion.
\begin{enumerate}
\item[(W3)]   The generating curve of $W$  has non-decreasing curvature (with respect to the inward pointing normal) as a function of arc length
on $\{z\ge 0\}$ as one moves in an upward direction.
\end{enumerate}
The condition (W3) will be referred to as the {\it curvature condition}.

Because of (W1), the Gauss map of $W$ defines a bijection of $W$ onto $S^2$.
Therefore quantities defined on $W$ can be expressed unambiguously on the sphere. In particular, the principal curvatures of $W$ with respect to the inward pointing normal, $\mu_i$, $i=1,2$  along the latitudinal and longitudinal directions, are given respectively by
$$1/\mu_2=\gamma-\nu_3\gamma'\:,\qquad 1/\mu_1=(1-\nu_3^2)\gamma''+1/\mu_2\:.$$
The uniform convexity means that $\mu_i>0$ holds. Moreover, the position vector $\xi_1$ on $W$ can be expressed in terms of the position vector  $\nu$ on $S^2$  as
\begin{equation}
\label{xi}
\xi_1 =\frac{1}{\mu_2} \nu +\gamma'(\nu_3)E_3\:.\end{equation}

Given an oriented embedded surface $X:\Sigma \rightarrow {\bf R}^3$ with unit normal field $\nu$, we define the {\it Cahn-Hoffman field} as the composition $\xi:=\xi_1\circ \nu$.   The first variation of energy defines the anisotropic mean curvature $\Lambda$ by
\begin{equation}
\label{deltaf}\delta {\cal F}=-\int_\Sigma \Lambda \delta X\cdot \nu \:d\Sigma +\oint_{\partial \Sigma} (\xi \times \delta X\cdot dX).
\end{equation}
Some local expressions for the anisotropic mean curvature are:
$$
\Lambda :=-{\rm trace}_\Sigma Ad\nu =-(\nabla \cdot D\gamma-2H\gamma), \quad A:=d\xi_1=(D^2\gamma+\gamma1)|_{S^2}\:.
$$
Here $D$, $D^2$ denote the gradient and Hessian operators acting on functions on $S^2$. These can always be identified with tensors on $\Sigma$ by parallel translation in ${\bf R}^3$ since $T_p\Sigma=T_{\nu(p)}S^2$. 
Another useful expression is
$$\Lambda=-\nabla \cdot \xi\:,$$
where $\nabla \cdot$ denotes the surface divergence.

The problem we consider here is to understand the evolution to equilibrium of a drop of material trapped between two horizontal
planes $\Pi_i$, $i=0,1$ located at heights $z=0$ and $z=h$ respectively. The surface energy is assumed to be of the form \rf{F}.
It will also be assumed that:
\begin{itemize}
\item the initial surface is axially symmetric,
\item the generating curve of the initial surface is a graph over the rotation axis.
\item throughout the evolution, the surfaces meets the planes $\Pi_i$ orthogonally.
\item the initial surface extends smoothly to an infinite periodic surface with period $h$.
\end{itemize}

Of course, the third assumption is a consequence of the fourth,
but we include it for clarity.

 It is easy to see that the axial symmetry is preserved by the evolution considered below. It is also consistent with the form of the known minimizer for sufficiently large volumes and neutral wetting along the interface of $\Sigma$ with the planes, \cite{KP2006}. For small volumes, the drop must disconnect during the minimization process. 

It will be shown below that if the initial volume is sufficiently large, then the generating curves of the evolving surfaces are also graphs for all times for which the flow is defined.

Although we  expect that a similar analysis can be carried out in arbitrary dimensions, we have decided to concentrate on the case of surfaces in ${\bf R}^3$ because of its obvious physical significance.

\section{Preliminary results}

We will restrict our attention to axially symmetric surfaces with vertical rotation axis. In addition, the generating curve will always be a graph of the form $r=r(z)$ and we write the surface as 
$X=(r(z)e^{i\theta}, z)$ where we have identified the first two coordinates on ${\bf R}^3$ with the complex plane ${\bf C}$. The outward pointing normal to the surface is $\nu=(1+r_z^2)^{-1/2}
(e^{i\theta}, -r_z)$.

Consider the flow:
\begin{equation}
\label{AMCF} \frac{\partial X}{\partial t} =(\Lambda-{\bar \Lambda})\nu\:, \end{equation}
with the boundary condition
  \begin{equation}
\label{AMCFb} \nu\cdot E_3 \equiv 0\:, \:{\rm on}\:\partial \Sigma\:. \end{equation}

 From \rf{amcf}, we obtain
\begin{equation}
\label{zt}
\partial_tz=(\Lambda-{\bar \Lambda})\:\nu_3\:.\end{equation}
It is evident from this formula that, in general, $z$ must depend on $t$, so, following \cite{E},  we write $z=z(\zeta, t)$ where
we choose $\zeta$ to lie in the interval $[0,h]$. We then express the immersion as
$$X=X(\zeta,t)=(r(z(\zeta,t),t), z(\zeta,t))\:.$$
We get
$$ r_z z_t+r_t=(\Lambda-{\bar \Lambda})\sqrt{1-\nu_3^2}\:$$
Combining this with \rf{zt}, we get, using that $\nu_3=-r_z/\sqrt{1+r_z^2}$,

$$r_t=(\Lambda-{\bar \Lambda})(\sqrt{1-\nu_3^2}-r_z\nu_3)
=(\Lambda-{\bar \Lambda})(\frac{1}{\sqrt{1+r_z^2}}+\frac{r_z^2}{\sqrt{1+r_z^2}})\:,$$
i.e.
\begin{equation}
\label{rt}
r_t=(\Lambda-{\bar \Lambda}) \sqrt{1+r_z^2}\:.
\end{equation}

The evolution for $r$ is identical with the evolution for $r$ under the flow
$$(\partial_t X)^\perp=(\Lambda-\Lambda_0)\nu\:,$$
where $\perp$ denotes the normal component.

Short time existence for \rf{rt} is standard, see for example Theorem 8.1.1 of \cite{Lu}. By the equation
\rf{lambda}, \rf{barl} given below, both $\Lambda$ and ${\bar \Lambda}$ are expressible in terms of $r$ and its derivatives.

  The  admissible variations for the variational problem are those which keep the surface between the planes, $\delta X\cdot E_3 \equiv 0$ on $\partial \Sigma$ , and fix the  three dimensional volume enclosed within the surface:
$$\int_\Sigma \delta X\cdot \nu\:d\Sigma=0\:.$$
The infinitesimal generator  $(\Lambda-{\bar \Lambda})\nu$
of the evolution \rf{AMCF}, \rf{AMCFb} clearly satisfies these conditions.

By \rf{AMCFb}, along $\partial \Sigma$, the normal $\nu$ is contained in the horizontal equator of $S^2$. By (W2), $\gamma'(\nu_3=0)=0$ holds and so by \rf{xi}, $\xi$ and $X_t$
are parallel along $\partial \Sigma$.  we then get from the first variation formula \rf{deltaf},
$$\partial_t {\cal F}[\Sigma_t]= -\int_\Sigma \Lambda (\Lambda-{\bar \Lambda}) \:d\Sigma
=-\int_\Sigma  (\Lambda-{\bar \Lambda})^2 \:d\Sigma \: ,$$
so the flow decreases the anisotropic energy.\\[2mm]
{\bf Remark} When the wetting is not neutral, it appears to be problematic to construct an analogous flow. In this case,
the boundary condition for the minimizer $\xi\cdot E_3\equiv c_i\ne 0$ on $\partial \Sigma \cap \Pi_i$, where $c_i$ are non zero constants related to the coupling constants for the wetting energy. This boundary condition is incompatible with the
flow given by \rf{AMCF} maintaining the drop between the planes  .

For the  boundary value problem we are considering, the morphology of the minimizer depends  the initial volume and whether or not the condition (W3) holds. If (W3) holds, it was shown in \cite{KP2006} that for  volumes greater than or equal to a critical value $V_0$, all stable equilibria must be cylinders. Below this value, any minimizer must either disconnect or loose contact with at least one of the supporting planes.

When (W3) does not hold, numerical simulations for a particular class of functionals, \cite{AKP},  show that anisotropic unduloids may occur as stable equilibria for a certain range of volumes. However, for large volumes, only cylinders occur and for sufficiently small volumes, there is no stable connected surface spanning the two supporting planes.
\\[2mm]

We will next show that of the initial volume is sufficiently large, the surfaces will not pinch off, i.e. the generating curves are bounded away from the rotation axis. 

The following lemma uses calibrations to show a minimizing property of graphs with $\Lambda=0$. It is well known and its proof is given for completeness.

\begin{lemma}
Let $\Sigma$ be a surface with zero anisotropic mean curvature which can be represented as a graph over a planar domain $\Omega$. Let $S$ be a piecewise smooth oriented surface which is contained in the cylinder $\Omega \times {\bf R}$ and which has the same boundary as $\Sigma$. Then,
$${\cal F}[\Sigma] \le {\cal F}[S]\:$$
holds.
\end{lemma}
{\bf Remark.} In saying that $\Sigma$ and $S$ share the same boundary, we mean that $\Sigma - S$  is the oriented boundary of an oriented 3-chain.\\
{\it Proof.} Let $\xi$ denote the Cahn-Hoffman field of $\Sigma$. The condition that $\Sigma$ has zero anisotropic mean curvature is expressed $\nabla \cdot \xi=0$.

Because $\Sigma$ is a graph, we can extend $\xi$ to a field ${\tilde \xi}$ on $\Omega \times {\bf R}$ by making ${\tilde \xi}$ constant on all vertical lines through points in $\Omega$. 
This field will satisfy  ${\tilde \nabla}\cdot {\tilde \xi}$, where ${\tilde \nabla}$ denotes the divergence operator on ${\bf R}^3$.  Then, by the Stokes' Theorem, we get
\begin{eqnarray*}
{\cal F}[\Sigma]&=&\int_\Sigma \gamma(\nu_\Sigma)\:d\Sigma=\int_\Sigma \xi \cdot \nu_\Sigma d\Sigma\\
&=&\int_\Sigma {\tilde \xi} \cdot \nu_\Sigma\: d\Sigma=\int_S {\tilde \xi}\cdot \nu_S \:dS\\
&\le& \int_S \gamma(\nu_S) \:dS={\cal F}[S]\:,
\end{eqnarray*}
where the inequality follows from \rf{W}.

\begin{cor} Let  ${\cal F}$ be an  axially symmetric anisotropic surface energy and let $C\subset \Pi_i$
be a circle.  Let $S$ be any piecewise smooth compact surface bounded by $C$ and which is contained in the cylinder over the disc bounded by $C$. Also, let $D$ be the flat disc bounded by $C$. Then
$${\cal F}[S]\ge {\cal F}[D]=\gamma(e_3)|D|\:$$
holds. \end{cor}
{\it Proof.} This follows immediately from the previous lemma using the fact that the disc has zero anisotropic mean curvature.

\begin{prop}
\label{P1} Let $\Sigma_0$ be an initial axially symmetric surface enclosing a volume $V$ and intersecting
the supporting planes orthogonally. Assume that 
\begin{equation}
\label{np}
{\cal F}[\Sigma_0] < \gamma(e_3)\frac{V[\Sigma_0]}{d}\: \end{equation}
holds and that the  flow \rf{AMCF} exists for all $t\in [0,T)$.
Then, then no pinching occurs. In particular
 \begin{equation}
\label{c0} r\ge c_0\:\end{equation}
holds for $t\in [0,T)$, where 
$$\gamma(e_3)\frac{V[\Sigma_0]}{d}-{\cal F}[\Sigma_0]=:\pi c_0^2\:.$$
\ \end{prop}
{\it Proof.} Recall that ${\cal F}[\Sigma_t]$ is non increasing. Let $C$ denote the cylinder between $\Pi_0$ and
$\Pi_1$ enclosing the same volume $V[\Sigma_0]$ as $\Sigma_0$. Let $r_C$ be the radius of this cylinder.  If $\rho:=r(\zeta, t_1)<c_0$ for some $t_1<T$ and some $\zeta$, we can find an annular part  of the surface $\Sigma_{t_1}$ bounded by two circles of radii $\rho$ and  radius $r_{t_1} >r_C$.  From this piece of $\Sigma_{t_1}$, we form a piecewise smooth disc type surface by filling the circle of radius $\rho$ with a disc. 

By the previous corollary, we obtain
$${\cal F}[\Sigma_0]+\pi c_0^2 > {\cal F}[\Sigma_{t_1}]+\pi\rho^2 > \gamma(e_3)\pi r_{t_1}^2\ge \gamma(e_3)\pi r_c^2  =\gamma(e_3)\frac{V(\Sigma_0)}{d}\:,$$
which gives a contradiction.\\[2mm]

\begin{lemma}
\label{barl} Assume the conditions (W1), (W2) hold
for $W$ and that $\Sigma$ is an axially symmetric surface intersecting the planes $\Pi_i$ orthogonally for which \rf{c0} holds. Then there exists $c_1=c_1(c_0)$ such that
$$ c_1\ge  |{\bar \Lambda }|$$
holds. If, in addition (W3) holds, then we have
$$0\ge {\bar \Lambda }\ge c_1>-\infty\:.$$

\end{lemma}
{\it Proof.}  The idea of the proof is essentially taken from \cite{Ath}.

We first note that the curvature condition (W3) can be expressed
\begin{equation}
\label{cc}\nu_3 \partial_{\nu_3}\mu_1 \ge 0 ,\:\forall \nu_3\:.\end{equation}

For an axially symmetric surface whose generating curve is a graph, the anisotropic mean curvature is given by

\begin{equation}
\label{lambda}\Lambda =\frac{k_1}{\mu_1}+\frac{k_2}{\mu_2}=\frac{r_{zz}}{\mu_1(1+r_z^2)^{3/2}}-\frac{1}{\mu_2 r(1+r_z^2)^{1/2}}\:.\end{equation}
and its average value is
\begin{equation}
\label{barl}
{\bar \Lambda} =\frac{\int_0^h\Lambda r(1+r_z^2)^{1/2}\:dz}{\int_0^h  r(1+r_z^2)^{1/2}\:dz}\:.\end{equation}

Recall that $\mu_i$ are the principal curvatures of the Wulff shape $W$ with respect to the inward pointing normal so $0<<\mu_i <\infty$ holds. It then follows easily that
$$0\ge \frac{\int_0^h \frac{k_2}{\mu_2}r(1+r_z^2)^{1/2}\:dz}{\int_0^h r(1+r_z^2)^{1/2}\:dz}
=\frac{\int_0^h \frac{-1}{\mu_2}\:dz}{\int_0^hr(1+r_z^2)^{1/2}\:dz} \ge c_2(c_0)\:.$$

From the boundary condition, we have $r_z=0$ on $\partial \Sigma_0$. Note that $r_{zz}/(1+r_z^2)
=(\arctan r_z)_z$. We get
\begin{eqnarray*}
\int_\Sigma \frac{k_1}{\mu_1} \:d\Sigma &=&2\pi \int_0^h \frac{rr_{zz}}{\mu_1(1+r_z^2)}\:dz\\
                                                             &=&2\pi \int_0^h \frac{r}{\mu_1}(\arctan r_z)_z\:dz\\
&=&-2\pi \int_0^h \frac{r_z}{\mu_1}(\arctan r_z)\:dz-2\pi \int_0^h \partial_z (\frac{1}{\mu_1}) r(\arctan r_z)\:dz\:.
\end{eqnarray*}
Note that, using Proposition \rf{P1}, we obtain
$$0\ge \frac{-2\pi \int_0^h \frac{r_z}{\mu_1}(\arctan r_z)\:dz}{2\pi\int_0^hr(1+r_z^2)^{1/2}\:dz}
\ge c_3(c_0)\:. $$
First assume that (W3) holds. Note that since $\nu_3=-r_z/\sqrt{1+r_z^2}$, we can write
\begin{eqnarray*}
-2\pi \int  \partial_z (\frac{1}{\mu_1}) r(\arctan r_z)\:dz&=& -2\pi \int \partial_{\nu_3} (\frac{1}{\mu_1}) r(\arctan \bigl(\frac{-\nu_3}{\sqrt{1-\nu_3^2}}\bigr))\:d\nu_3\\
&=& -2\pi \int  (\frac{\partial_{\nu_3} \mu_1}{\mu^2_1}) r(\arctan \bigl(\frac{\nu_3}{\sqrt{1-\nu_3^2}}\bigr))\:d\nu_3\\
&\le& 0\:
\end{eqnarray*}
by \rf{cc}, since $\nu_3$ and $\arctan (\nu_3 (1-\nu_3^2)^{-1/2})$ have the same sign. Also, since $\partial_{\nu_3} \mu_1/\mu_1^2$ is uniformly bounded, we obtain that
$$0\ge \frac{\int_\Sigma \frac{k_1}{\mu_1} \:d\Sigma}{\int_\Sigma \:d\Sigma} \ge c_4(c_0)\:,$$
and the result follows by combining this with the previous inequalities. 

It (W3) is not assumed to hold, then we easily obtain $|(k_1/\mu_1)|\:d\Sigma \le c r \:dz\:d\theta$ for a constant $c$ while the integrand in the denominator is counded below by $r$.  {\bf q.e.d}
\\[2mm]

\section{Evolution equations}
Again, the governing evolution equation is
\begin{equation}
\label{amcf}
\partial_tX=(\Lambda-{\bar \Lambda})\:\nu\:.\end{equation}

For any smooth variation $\delta X=\psi \nu +T$ of a surface $X$, the corresponding pointwise variation of the anisotropic mean curvature is
\begin{equation}
\label{J0} \delta \Lambda=J[\psi]+\nabla \Lambda\cdot T\:.\end{equation}
Here $J$ is the self=adjoint elliptic operator
given by
\begin{equation}
\label{J}
J[u]=\nabla\cdot A\nabla u+\langle Ad\nu,d\nu \rangle u\:.\end{equation}
Using that the evolution of $X$ is given by
\rf{AMCF}, we obtain
\begin{equation}
\label{PL}
\Lambda_t-\nabla\cdot A\nabla \Lambda=\langle d\nu,Ad\nu \rangle(\Lambda-{\bar \Lambda})\:.
\end{equation}

  We recall from
\cite{KP2006} that the normal $\nu$ satisfies the equation
\begin{equation}
\label{jnu}
\nabla\cdot A\nabla \nu_j+\langle Ad\nu,d\nu \rangle \nu_j=-\nabla \Lambda \cdot E_j\:,\:j=1,2,3\:.\end{equation}
This is a consequence of \rf{J0} and the translation invariance of the functional.

To compute the evolution of the normal, we use that for $\delta X=\psi\nu +T$,
one has $\delta \nu=-\nabla \psi+d\nu T$.
Since $X_t=(\Lambda-{\bar \Lambda})\nu$, we get
\begin{equation}
\label{nut}\nu_t=-\nabla \Lambda\end{equation}

We define the parabolic operator
$$P[f]=f_t-\nabla\cdot A\nabla f\:.$$
\begin{lemma}
Define $\omega:=\sqrt{1+r_z^2}=1/\sqrt{1-\nu_3^2}\:$. Then

\begin{equation}
\label{pomega}
P[\omega]=\frac{-2\nabla \omega\cdot A\nabla \omega}{\omega}+\frac{\omega}{r^2\mu_2}-\langle A\:d\nu,d\nu \rangle \omega\:.\end{equation}
\end{lemma}
{\it Proof.} By combining \rf{jnu} and \rf{nut},
we obtain the vector equation
$$P[\nu]=-\nabla \Lambda\:.$$
For the surface of revolution, we have
$$\nu=(\sqrt{1-\nu_3^2}\:e^{i\theta}, \nu_3)\:,$$ where we have identified the space of the first two coordinates with the complex plane.
It follows from \rf{jnu} that
\begin{equation}
\label{eq}\nabla \cdot A\nabla (\sqrt{1-\nu_3^2}\,)e^{i\theta})=-\langle d\nu , Ad\nu\rangle (\sqrt{1-\nu_3^2}e^{i\theta}-\underline{\nabla \Lambda}\:,\end{equation}
where the last term denotes the projection of $\nabla \Lambda$ onto ${\bf R}^2$.

The metric on $\Sigma$ is $dS^2=(1+r_z^2)dz^2+r^2d\theta^2:=\alpha_1\otimes
\alpha_1+\alpha_2\otimes
\alpha_2$  If $f(z)$ is a differentiable function, then $\nabla f\cdot \nabla \theta d\Sigma =df\wedge *d\theta=f_zdz\wedge *\alpha_2/r=f_zdz\wedge (-\omega /r)dz=0$, and therefore $\nabla f\cdot \nabla \theta =0$. Also $\Delta \theta\:dS =d*d\theta
=d *(\alpha_2/r)=d(-\omega /r)dz =0$. Finally $|\nabla \theta |^2dS=d\theta \wedge *d \theta=d\theta \wedge (-\omega/r)dz=(1/r^2)\:d\Sigma$ so $|\nabla \theta|^2=1/r^2$.

Using these formulas to expand out the left hand side of \rf{eq} and using that $\nu_3$ only depends on $z$,  we have
\[
\begin{array}{l}
\disp \nabla \cdot A\nabla  (\sqrt{1-\nu_3^2}e^{i\theta}) \\
\disp \qquad = e^{i\theta}\nabla \cdot A\nabla  (\sqrt{1-\nu_3^2}) +ie^{i\theta}\frac{\sqrt{1-\nu_3^2}}{\mu_2}(i|\nabla \theta|^2+\Delta \theta)\\
\disp \qquad = e^{i\theta}\nabla \cdot A\nabla  (\sqrt{1-\nu_3^2}) -e^{i\theta}\frac{\sqrt{1-\nu_3^2}}{r^2 \mu_2}\:.
\end{array}
\]
Combining this with \rf{eq}, we get
\begin{equation}
\label{eq2}
\nabla\cdot A\nabla (\sqrt{1-\nu_3^2})=\frac{\sqrt{1-\nu_3^2}}{\mu_2r^2}-\langle Ad\nu, d\nu \rangle\sqrt{1-\nu_3^2}-e^{-i\theta} \underline{\nabla \Lambda}\:.
\end{equation}
From this, it is easy to obtain
\[
\begin{array}{l}
\disp \nabla \cdot A\nabla \omega \vphantom{\frac{1}{2}}\\
\disp \qquad =\nabla \cdot A\nabla \frac{1}{\sqrt{1-\nu_3^2}}
-(1-\nu_3^2)^{-1}\nabla\cdot A\nabla (\sqrt{1-\nu_3^2}) \\
\disp \qquad \quad  -\nabla (1-\nu_3^2)^{-1}\cdot A\nabla \sqrt{1-\nu_3^2} \vphantom{\frac{1}{2}} \\
\disp \qquad = -\frac{\omega}{\mu_2r^2}+ \langle Ad\nu, d\nu \rangle\omega+
e^{-i\theta}\omega^2\underline{\nabla \Lambda}-\nabla (1-\nu_3^2)^{-1}\cdot A\nabla \sqrt{1-\nu_3^2}\\
\disp \qquad = -\frac{\omega}{\mu_2r^2}+ \langle Ad\nu, d\nu \rangle\omega+
e^{-i\theta}\omega^2\underline{\nabla \Lambda}+\frac{2}{\omega}\nabla \omega \cdot A\nabla \omega\\
\disp \qquad = (*) \vphantom{\frac{1}{2}} \:.
\end{array}
\]
From \rf{nut}, we get $ (\sqrt{1-\nu_3^2})_t=-e^{-i\theta}\underline{\nabla \Lambda}$ which implies
$$\omega_t=(1-\nu_3^2)^{-1} e^{-i\theta}\underline{\nabla \Lambda}=\omega^2
e^{-i\theta}\underline{\nabla \Lambda}\:.$$
Combining this with (*) yields \rf{pomega}
\bigskip

We state the following well known Maximum Principle, \cite{PW}.

\begin{theorem}
If the operator $u\mapsto \nabla \cdot A\nabla[u]$ is elliptic and $h:M\times [t_1,t_2]\rightarrow {\bf R}$ satisfies
$$(\partial_t-\nabla \cdot A\nabla)h\le 0\:,$$
holds, then
$$\max_M h(\cdot,t)\le \max_M h(\cdot ,t_1)\:.$$
\end{theorem}
In the case we are considering, ellipticity of $u\mapsto \nabla \cdot A\nabla[u]$ follows from the convexity condition W1.
\begin{prop}
\label{oprop} Assume that the initial surface satisfies the estimate
\rf{np}. Then the generating curve is a graph for all time $t$ for
which the evolution \rf{AMCF} exists. \end{prop}
{\it Proof.}  As in \cite{Ath}, we use the boundary condition \rf{AMCFb} and the assumptions on the  initial curve to extend the surfaces generated by the flow to periodic surfaces.

For  any constant $c$, we obtain from \rf{pomega}
\begin{equation}
(\log (\omega)-ct)_t-\nabla\cdot A\nabla(\log (\omega)-ct)=\frac{1}{r^2\mu_2}-c- \frac{\nabla \omega\cdot A\nabla \omega}{\omega^2}-\langle A\:d\nu,d\nu\rangle\:.\end{equation}
By Proposition \rf{P1}, there is a constant $c_0=c_0(t_1)$ such that $r\ge c_0$ holds for $t\in [t,t_1]$. Hence, for a suitable constant $c$, we have $(\log (\omega)-ct)_t-\nabla\cdot A\nabla(\log (\omega)-ct)\le 0$.
It follows from the Maximum Principal, that
\begin{equation}
\label{max} \max_{0\le z\le d, 0\le t\le t_1} (\log (\omega)-ct) \le \max_{0\le z\le d, t=0}\log \omega\:.
\end{equation}

 From the definition of $\omega$ we have from \rf{max} that there exist constants $c_2(\Sigma_0)$, $a_2(\Sigma_0)$ such  that if
the flow exists for $t\in [0,T)$, then
\begin{equation}
\label{wbound}
\omega \le c_2e^{a_2t}:\,\end{equation}
for $0\le t< T$.

Recalling the definition of $\omega$, the result follows. {\bf q.e.d.} 
\\[4mm]
\begin{prop}
Assume the conditions (W1) and (W2) hold. Assume that the flow \rf{AMCF} is defined for $0\le t < T$. Then the curvatures of the surfaces $X(t)$ remain bounded, i.e. there exists a constant $c_3(T)$.
with
\begin{equation}
\label{dn}
 k_1^2(t)+k_2^2(t)=:||d\nu||^2(t)\le c_3(T)\:,\forall  t, 0\le t< T\end{equation}
where $k_i(t)$ denote the principal curvatures at time $t$.
\end{prop}
First note that $k_2=-(r\omega)^{-1}$. By Proposition \rf{P1} we have $r \ge c_0(T)$ and $c_4(T)\ge \omega \ge 1$ by the previous
proposition and the definition of $\omega$. It follows that $|k_2|\le c_5(T)$. Since $\mu_i$, $i=1,2$, are uniformly bounded below and above,
($0<a\le \mu_i\le 1/a$ for some $a\in {\bf R}$), it is enough to show the existence of a bound
\begin{equation}
\label{boundl}
\Lambda^2 \le c_6(T)\end{equation}
for $0\le t <T$ and then \rf{dn} will follow.

We recall the standard formula:
\begin{equation}
\label{pproduct}
P[fg]=fP[g]+gP[f]-2\nabla f\cdot A\nabla g\:,\end{equation}
for sufficiently smooth functions $f$ and $g$.
The endomorphism field  $A=D^2\gamma+\gamma {\rm id}$ is self-adjoint so that the last term is symmetric in $f,g$.

From the last equation \rf{PL} and \rf{pomega}
we obtain
\begin{equation}
\label{pw}
P[\omega^2]=\frac{2\omega^2}{\mu_2r^2}-2\langle Ad\nu, d\nu\rangle \omega^2-6\nabla \omega A\nabla \omega\:,
\end{equation}
\begin{equation}
\label{pL}P[\Lambda^2]=2\langle Ad\nu, d\nu\rangle \Lambda (\Lambda-{\bar \Lambda})-2\nabla \Lambda\cdot A\nabla \Lambda \:.\end{equation}

For a suitable $C^2$ function $h(x)$ and $B\in {\bf R} $ both of  which we determine later, we get
\[
\begin{array}{l}
\disp P[h(\omega^2)(\Lambda^2+B)] \vphantom{\dfrac{1}{2}}\\
\disp \qquad =(\Lambda^2+B) h'(\omega^2)\bigl(\frac{2\omega^2}{\mu_2r^2}-2\langle Ad\nu, d\nu\rangle \omega^2-6\nabla \omega A\nabla \omega \bigr)\\
\disp \qquad \quad -\;(\Lambda^2+B)h''(\omega^2) \nabla \omega^2\cdot A\nabla \omega^2 \vphantom{\dfrac{1}{2}}\\
\disp \qquad \quad  + \;h(\omega^2)\bigl(2\langle Ad\nu, d\nu\rangle \Lambda (\Lambda-{\bar \Lambda})-2\nabla \Lambda\cdot A\nabla \Lambda\bigr)
-2h'(\omega^2)\nabla \omega^2\cdot A\nabla \Lambda^2. \vphantom{\dfrac{1}{2}}
\end{array}
\]

We take $h(x):=e^{ax}$ for a constant $a$ to be determined later, and get
\[
\begin{array}{l}
 \disp P[e^{a\omega^2}(\Lambda^2+B)] \vphantom{\dfrac{1}{2}}\\
 \disp \qquad = e^{a\omega^2}\biggl(
(\Lambda^2+B) a\bigl(\frac{2\omega^2}{\mu_2r^2}-2\langle Ad\nu, d\nu\rangle \omega^2-6\nabla \omega A\nabla \omega \bigr) \\

\disp \qquad \quad -(\Lambda^2+B)a^2 \nabla \omega^2\cdot A\nabla \omega^2 
  + 2\langle Ad\nu, d\nu\rangle \Lambda (\Lambda-{\bar \Lambda}) \vphantom{\dfrac{1}{2}} \\
  
\disp \qquad \quad -2\nabla \Lambda\cdot A\nabla \Lambda
-2a\nabla \omega^2\cdot A\nabla \Lambda^2\biggr) \vphantom{\dfrac{1}{2}} \\
\disp \qquad = e^{a\omega^2}\biggl( \frac{2a(\Lambda^2+B)}{\mu_2r^2}
+\bigl[ -2a(\Lambda^2+B)+2\Lambda^2-2\Lambda {\bar \Lambda}\bigr]\langle Ad\nu, d\nu\rangle\\
\disp \qquad \quad +\bigl[-6a^2(\Lambda^2+B)\nabla \omega^2 \cdot A\nabla \omega^2
-2\nabla \Lambda\cdot A\nabla \Lambda -2a\nabla \omega^2\cdot  A\nabla \Lambda^2\bigr]             \biggr) \vphantom{\dfrac{1}{2}} \\
\disp \qquad = (**) .
\end{array}
\]

We have
$$-2a(\Lambda^2+B)+2\Lambda^2-2\Lambda {\bar \Lambda} \le
-2a(\Lambda^2+B)+3\Lambda^2+{\bar \Lambda}^2\:,$$
using $2xy\le x^2+y^2$. Because of Lemma \rf{barl}, we have the existence of a constant $c_1(T)$ with ${\bar \Lambda}^2 \le (c_1(T))^2$.
Hence by choosing $a,B>>0$, we get that the term in (**) which includes the factor $\langle Ad\nu, d\nu \rangle $ is non positive.

Next notice that since the tensor $A$ is positive definite and self adjoint,
at each fixed point $p\in\Sigma$, we have an inner product defined by
$(u,v)=Au\cdot v$, $u,v\in T_p\Sigma$.  We then have using $2|(u,v)|\le
(u,u)+(v,v)$,
\begin{eqnarray*}
2a|\nabla \omega^2\cdot A\nabla \Lambda^2|&=&
4a|\Lambda||\nabla \omega^2\cdot A\nabla \Lambda|\\
&=& 4a|\Lambda||(\nabla \omega^2,\nabla \Lambda)\\
&\le& 4a^2\Lambda^2 (\nabla \omega^2,\nabla \omega^2)+(\nabla \Lambda,\nabla \Lambda)\\
&=& 4a^2 \Lambda^2 \nabla \omega^2 \cdot A\nabla \omega^2+
\nabla \Lambda\cdot A\nabla \Lambda\:.
\end{eqnarray*}
This means that both terms in (**) between the square brackets are non positive for suitable choices of the constants $a$ and $B$ and so we can conclude that
$$P[e^{a\omega^2}(\Lambda^2+B)]\le e^{a\omega^2} \frac{2a(\Lambda^2+B)}{\mu_2r^2}\:,$$
holds. Now recall that by Proposition \rf{P1}, we have $r\ge c_0(\Sigma_0)>0$.
It then follows that
\begin{eqnarray*}
P[e^{-Mt}e^{a\omega^2}(\Lambda^2+B)]&=& e^{-Mt}\biggl(
-Me^{a\omega^2}(\Lambda^2+B)+P[e^{a\omega^2}(\Lambda^2+B)]\biggr)\\
&\le&e^{-Mt}\biggl(
-Me^{a\omega^2}(\Lambda^2+B)+e^{a\omega^2} \frac{2a(\Lambda^2+B)}{\mu_2r^2}\biggr)\\
&\le& 0\:,
\end{eqnarray*}

for a sufficiently large constant $M$. Recalling from Proposition \rf{oprop}
that $\omega \le \sqrt{1+ c_2^2(T)}$, we find that \rf{boundl} follows from the Maximum Principle. {\bf q.e.d.}

\section{Higher order regularity}

In this section we obtain bounds on higher order derivatives of the surface. In similar problems involving mean curvature flow, these are obtained from bounds on higher order derivatives of the second fundamental form $|d\nu|^2$. In our case, a nice evolution equation for this quantity is unavailable. However, since we are only working with  surfaces of revolution, we can use \rf{pL} instead  to obtain the desired bounds.

\begin{prop}
\label{pbound}
For positive constants $\sigma_m=\sigma_m(T)$, $\tau_m=\tau_m(T)$, there holds
\begin{equation}
\label{pb1}
P[|\nabla ^m\Lambda|^2] \le -\sigma_m|\nabla ^{m+1}\Lambda|^2+\tau_m(1+|\nabla ^m\Lambda|^2)\:,
 \tau_m=\tau_m(|\nabla ^j \Lambda|^2)\:, j < m\end{equation}
 \end{prop}
 for $0\le t \le T$.
\begin{lemma}
\label{pb3}
Let $f=f(s)$ be a sufficiently smooth function. Then
\begin{equation}
\label{pb4}
|\nabla ^m f(s)\le |\nabla_s ^m f(s)|^2+ E_1\:,\end{equation}
where $E_1$ only depends on $|\nabla^j f|$, $ j\le m-1$.
\end{lemma}
{\it Proof.} The proof is straightforward and is left to the reader.

{\it Proof of Proposition \rf{pbound}}
At a fixed time $t$, we let $\omega_1:=ds$, $\omega_2=rd\theta$ and we write, using the summation convention, 
$$d\Lambda=\Lambda_j \omega_j\:,\:\:\nabla d\Lambda =\Lambda_{jk}\omega_j\omega_k\:,...$$
Then at any time $t$, we can write
 \begin{equation}
 \label{Lm}|\nabla^m \Lambda|^2 = |\Lambda_{j_1....j_m}|^2 \prod_{k=1..m} g^{j_k j_k}\:,\end{equation}
 where, although the coframe is fixed, all quantities depend on $t$.  Throughout the evolution,
 all the surfaces are axially symmetric, so $g^{ii}=1/g_{ii}$.  For a normal variation $X(t)=X+t\psi \nu+...$, the first variation of the metric is $\partial_t g_{ij}=-2\psi h_{ij}$ where $h_{ij}$ are
 the coefficients of the second fundamental form.  Since all the surfaces are axially symmetric, $g_{ij}=\delta_{ij}g_{ii}$ and $h_{ij}=\delta_{ij}h_{ii}$ for all $t$. Using this, we obtain from \rf{Lm}
 \begin{eqnarray*}
 \partial_t |\nabla ^m\Lambda|^2 &=&4H(\Lambda-{\bar \Lambda})|\nabla ^m\Lambda|^2+
 2\Lambda_{j_1....j_m}( \partial_t\Lambda_{j_1....j_m}) \prod_{k=1..m} g^{j_k j_k}\\
 &=&4H(\Lambda-{\bar \Lambda})|\nabla ^m\Lambda|^2+2\langle \nabla^m \Lambda, \nabla^m
 J[\Lambda] \rangle\\
 &\sim& 2\langle \nabla^m \Lambda, \nabla^m
(\nabla \cdot A\nabla \Lambda) \rangle\:,
 \end{eqnarray*}
 where $\sim$ means that the quantities are equal up to terms of orders less than or equal to $m$.
 In the given frame, for a function $f(s)$, 
\begin{equation}
 \label{div}
\nabla \cdot A\nabla f =\frac{1}{r\omega} \partial_s(\frac{rf_s}{\mu_1})\:.\end{equation}
 In particular, since $\Lambda$ only depends on $s$, we get
 \begin{equation}
\label{pb5}
 \partial_t |\nabla ^m\Lambda|^2 \sim \frac{2}{\omega \mu_1} (\nabla^m_s \Lambda) (\nabla_s^{m+2} \Lambda)+\Phi_m(\nabla^m_s \Lambda)(\nabla_s^{m+1} \Lambda)\:,
 \end{equation}
 for suitable functions $\Phi_m$ which are bounded by lower order derivatives.
 
Since $|\nabla^m \Lambda|^2$ only depends on $s$, applying \rf{div} gives
\begin{eqnarray*}
\nabla \cdot A\nabla |\nabla^m \Lambda|^2& =&\frac{1}{r\omega} \partial_s(\frac{r\partial_s|\nabla^m \Lambda|^2}{\mu_1})\\
&=& \frac{1}{r\omega} \partial_s(\frac{2r}{\mu_1}\langle \nabla^m\Lambda, \nabla_s\nabla^m \Lambda \rangle)\\
&=& \frac{2}{\omega \mu_1} \langle \nabla^m \Lambda, \nabla_s \nabla_s \nabla^m \Lambda\rangle +\frac{2}{\omega \mu_1}|\nabla_s^{m+1} \Lambda|^2 \\
 && \qquad +\frac{1}{r\omega}(\frac{2r}{\mu_1})_s \langle \nabla^m \Lambda,  \nabla_s \nabla^m \Lambda\rangle\:.
\end{eqnarray*}
From the last equation, \rf{pb5} and Lemma \ref{pb3}, we can write
$$ P[\:|\nabla ^m\Lambda|^2\:]\sim -\frac{2}{\omega \mu_1}|\nabla_s^{m+1} \Lambda|^2
+{\tilde \Phi}_m  \langle \nabla^m \Lambda,  \nabla_s \nabla^m \Lambda\rangle\:.$$
for a suitable function ${\tilde \Phi_m}$ which depends only on lower order derivatives.
 The second term can be bounded 
 $$|{\tilde \Phi}_m  \langle \nabla^m \Lambda,  \nabla_s \nabla^m \Lambda\rangle|
 \le \frac{\epsilon^2}{2}|\nabla_s^{m+1} \Lambda|^2+ \frac{\epsilon^{-2}}{2}|{\tilde \Phi}_m||^2\nabla^m\Lambda|^2
 \sim \frac{\epsilon^2}{2}|\nabla_s^{m+1} \Lambda|^2\:.$$
 By applying lemma \rf{pb3}, we obtain
$$ P[\:|\nabla ^m\Lambda|^2\:]\sim -\sigma_m|\nabla_s^{m+1} \Lambda|^2\:,$$
where for suitable $\epsilon$, we can take, using \rf{wbound}, 
\begin{equation}
\label{sbound}
 \sigma_m:=\big(\max(\mu_1)\sqrt{1+c_2^2(T)}\:\:\bigr)^{-1}\:.\end{equation} This shows that \rf{pb1} holds. {\bf q.e.d.}\\[4mm]

\begin{theorem}
\label{lmbound}

There exist constants $C_m(T)$ such that
\begin{equation}
\label{blm}
|\nabla^m \Lambda |^2 \le C_m(T)\:, \forall t < T\:.\end{equation}
\end{theorem}
 {\it Proof of Th. \rf{lmbound}.}  Denote by $S(m)$  the statement
 \rf{blm}. The statement $S(0)$ is just \rf{boundl} which was shown
 above. We assume $S(j)$ holds for all $j<m$.

  Let $\Psi_m:=|\nabla^m \Lambda |^2$. By \rf{pproduct}, \rf{pw} and \rf{pb1}, we find, for positive  constants $a$ and $B$
  \[
  \begin{array}{l}
 \disp P[e^{a\omega^2} (\Psi_m+B)] \vphantom{\frac{1}{2}}\\
 \disp \qquad \le  e^{a\omega^2}\biggl(  
 -\sigma_m\Psi_{m+1}+\tau_m(1+\Psi_m)\\
 \disp \qquad \quad +(\Psi_m+B)( a\bigl[\frac{2\omega^2}{\mu_2r^2}-2\langle Ad\nu, d\nu\rangle \omega^2-6\nabla \omega A\nabla \omega \bigr]-a^2 \nabla \omega^2 \cdot A\nabla \omega^2)\\
 \disp \qquad \quad - 2a \nabla \Psi_m\cdot A\nabla \omega^2
 \biggr)\:.
 \end{array}
 \]
 We can estimate the last term above using
 $$2a|\nabla \Psi_m\cdot A\nabla \omega^2| \le
 \epsilon^2\Psi_{m+1} +\epsilon^{-2}a^2 \nabla \omega^2 \cdot A\nabla \omega^2 \:.$$
 
 We choose $\epsilon <\sigma_m(T)$ and choose $B>>0$ so that
 $-6B\nabla \omega A\nabla \omega+\epsilon^{-2}a^2 \nabla \omega^2 \cdot A\nabla \omega^2<0$ holds.  This is possible since $\omega$ is bounded by \rf{wbound}. Recalling that $(\omega /r)$ is also bounded in any finite time interval, so we arrive at
$$ P[e^{a\omega^2} (\Psi_m+B)] < C e^{a\omega^2} (\Psi_m+B)\:,$$ for $t\le T$ and a suitable constant $C=C(T)$. Here we are using the induction hypothesis since the constant $\tau_m=
\tau_m(|\nabla^j\Lambda|^2)$, $j<m$, have been absorbed into the costant $C(T)$. It then follows that
$$P[e^{-Ct}e^{a\omega^2} (\Psi_m+B)]<0\:$$
holds for $0\le t\le T$, so the result follows from the Maximum Principle. {\bf q.e.d.}

 \begin{prop} If \rf{blm} holds then for suitable constants, the principal curvatures satisfy 
 \begin{equation}
 \label{dnm}
 |\nabla ^m k_1|\le c_m\:, \quad |\nabla^m
 k_2|\le c_m'\end{equation} hold.
 \end{prop}
 {\it Proof.} For a surface of revolution, the Codazzi equations reduce to 
 $$ (k_2)_s=(k_1-k_2)r_s/r\:.$$
 It follows by an easy induction argument that an upper bound for 
 $|\nabla^m_s k_2|$ and hence $|\nabla^m k_2|$,  can be obtained from upper bounds on 
 $|\nabla^j k_i|$, $i=1,2$, $j<m$. For this one needs \rf{c0} and the fact that derivatives of $r$ of order $j$ have upper bounds which depend on derivatives of $k_1$ of order $\le j-2$ since $r_{ss}=k_1z_s$.
 
    For a surface of revolution, the {\it anisotropic principal curvatures} are the functions $\lambda_j:=k_j/\mu_j$, $j=1,2$. 
    These are the eigenvalues of the differential $d\xi:T\Sigma \rightarrow TW$. Note that $\Lambda =\lambda_1+\lambda_2$. 
    
    We have $(k_j)_s=(\lambda_j)_s \mu_j +\lambda_j\mu_j'(\nu)_3(-k_1)
    z_s$ from which it follows easily that
   $$|\nabla ^m_s k_j|\le (\max \mu_j) | \nabla^m_s \lambda_j| +E_2\:,$$
   $$|\nabla ^m_s \lambda_j|\le (\max (1/\mu_j)) | \nabla^m_s k_j| +E'_2\:,$$
  where $E_2, E_2'$ depends on derivatives of the $k_j$'s of order
  less than or equal to $m-1$.
  
  Finally, we have from the definition of $\Lambda$,
 \begin{eqnarray*} |\nabla^m \lambda_1|&\le& |\nabla^m \Lambda |+|\nabla ^m \lambda_2|\\
 &\le&|\nabla^m \Lambda|+ (\max (1/\mu_j)) | \nabla^m_s k_2| +E'_2\\
&\le&|\nabla^m \Lambda|+E_3+E_2'\:, 
 \end{eqnarray*}
where $E_3$ (and $E_2'$) only depends on $|\nabla^{j}k_i|$, $j\le m-1$. Using induction, \rf{blm}  and \rf{dn}, the result follows. {\bf q.e.d.}\\[4mm]
\begin{theorem}
\label{lt}
Assume that the initial surface satisfies \rf{np}. Then the flow exists
for all time.
\end{theorem}
{\it Proof.} The result follows by a standard argument. Assume 
the flow exists on a finite time interval $[0, T_{\rm max})$. Because of the uniform estimates given by \rf{dn} and \rf{dnm}, the flow can also be extended smoothly to $t=T_{\rm max}$. Then, by applying the local existence result, the flow can be extended to $[0, T_{\rm max}+\epsilon)$ for some 
$\epsilon>0$. {\bf q.e.d.}\\[4mm]

{\bf Remark} We show $L^2$ convergence of the anisotropic mean curvature to its mean.
Note that 
$$\partial_t {\cal F}[\Sigma_t]= -\int_\Sigma (\Lambda-{\bar \Lambda})^2\:d\Sigma\:,$$
and therefore 
$${\cal F}[\Sigma_0]\ge \int_0^\infty \int_{\Sigma_t}  (\Lambda-{\bar \Lambda})^2\:d\Sigma_t\:dt\:.$$
In particular
$$ \int_\Sigma (\Lambda-{\bar \Lambda})^2\:d\Sigma_t \rightarrow 0\:,$$
as $t\rightarrow \infty$.

\section{Numerical Results}

Based on the different descriptions of the evolving surface, different methods can be used to numerically solve a surface evolution equation. These include parametric, level set, and phase field methods (see \cite{DDE05}). Each method has its own advantage and disadvantage. We choose the parametric method, which basically is a front tracking method, namely, the surface is evolved and tracked according to the surface evolution equation \eqref{rt}. 

Instead of using fully implicit schemes for the temporal discretization of equation \eqref{rt}, a semi--implicit backward Euler method is used.
The key of time forwarding in the semi--implicit scheme is to approximate the nonlinear terms in the equation  by using the previously computed approximated solution, while the linear terms still need to be solved implicitly. Therefore we can avoid solving systems of nonlinear equations (for example, using the Newton's method) at each time step and thus the computational cost can be reduced. On the other hand, the implicit feature will increase the stability of the scheme so the restrictions on time step sizes can be loosened.

As for the spatial discretization, a second--order finite difference method can be used. Finite element method using piecewise linear functions can also be used, and if the mass is lumped, it will be equivalent to the finite difference method. However, based on our numerical experiments with the above two methods, we choose to present the method of using cubic spline approximations.  The spline approximation not only provides a higher order method, it also ensures the continuity of the second order derivatives across the spatial nodes, which we think it is important in the approximation of curvature flows.
On the other hand, since the curvature flow we study here is essentially one dimensional, as we will see below, when using the spline approximation, the resulting linear system that we need to solve is very sparse and the computational cost will be of the same order as the cost using finite difference or finite element method.

\subsection{Approximation by Splines}

In this section, we discuss numerical methods for simulating the initial boundary value problem  \rf{AMCF} and \rf{AMCFb}.  For a representative class of examples, we concentrate on the Rapini--Papoular functionals given by $\gamma =1+\epsilon \nu_3^2$, where $-1\le \epsilon \le 1$. 

We use the clamped cubic spline approximation to numerically solve the evolution equation.
Let $T>0$ and let $N, K$ be positive integers. Let $0=z_1\le\cdots\le z_N\le z_{N+1}=1$ be a partition of the interval $[0,1]$ and let $0=t_0\le\cdots\le t_K\le =T$ be an equally spaced partition of the time interval $[0,T]$. 

We use the following notation
\[ 
   \disp r_n^{k}=r(t_k,z_n), \quad (r_z)_n^k = r_z(t_k,z_n), \quad 
   (r_{zz})_n^k = r_{zz}(t_k,z_n), \quad
   \disp \overline{\Lambda}_n^k = \overline{\Lambda}(t_k, z_n),\\
\]
where $0\le k\le K$ and $1\le n\le N+1$. We also set 
\[ \disp \tau = T/K , \qquad Q(r(t,z)) = 1 + (r_z(t,z))^2\,. \]

Using the semi--implicit backward Euler method and the above notations, the discretized evolution equation is
\be \label{devolution}
\disp \dfrac{r_n^{k+1}- r_n^{k}}{\tau}  = \Big(\dfrac{(r_{zz})_n^{k+1}}{\mu_1 [Q((r_z)_n^k)]^{3/2}} - \dfrac{1}{\mu_2 [Q((r_z)_n^k)]^{1/2}} - \overline{\Lambda}_n^{\,k} \Big) \sqrt{Q((r_z)_n^k)}\,, 
\ee
for $0\le k\le K$ and $1\le n\le N+1$. 

For an approximation of the generating curve of the axially surface $r(t,z)$, we seek a spline function $S(t,z)$  that satisfies equation \eqref{devolution} and the following properties :
\bi 
  \item[(i)] $S_n(t,z)= S(t,z)\big\lvert_{[z_n,z_{n+1}]}$, the restriction of $S(t,z)$ on interval $[z_n,z_{n+1}]$, $\, 1\le n\le N$, is a polynomial of degree no more than $3$.
  \item[(ii)] The first and second order partial derivatives of $S$ with respect to $z$ exist at nodes $z_1, \cdots, z_{N+1}$ and are continuous at the internal nodes
 $z_2, \cdots, z_N$.
   \item[(iii)] $\disp \dfrac{\partial S}{\partial z}(t,z_1)=\alpha=\dfrac{\partial r}{\partial z}(t,0)$ and $\disp \dfrac{\partial S}{\partial z}(t,z_{N+1})=\beta= \dfrac{\partial r}{\partial z}(t,1)$, where $\alpha=\beta=0$ if the contact angles of the surface at the top and bottom planes are required to be right angles.
\ei

At time $t=t_{k+1}$, $\,0\le k\le K-1$, the unknown for equation \eqref{devolution} is denoted by a column vector of length $2N$:
\be\label{unknown} \disp \vec{x}^{\,k+1}=(r_1^{k+1}, \cdots, r_{N+1}^{k+1}; d_2^{k+1}, \cdots, d_N^{k+1})^T \, 
\ee
where 
\[  \disp r_n^{k+1}=S(t_{k+1},z_n), \quad d_n^{k+1} =  \dfrac{\partial S}{\partial z}(t_{k+1},z_n),\quad 1\le n \le N+1\,.\]

To derive a linear system of $\vec{x}^{\,k+1}$ from equation \eqref{devolution}, for each $
1\le n \le N$, we adopt the following notations:
 \[ \ba{lll}
 \disp S_n^{k}(z) &= & \disp S_n(t_k, z) = S(t_k,z)\Big\lvert_{[z_n,z_{n+1}]}(t_k,z),
 \quad 0\le k \le K,\\
 \disp h_n & =  \disp & z_{n+1}-z_{n},\vphantom{\dfrac{1}{2}} \\
 \disp \delta_n & = & \disp \dfrac{r_n^{k+1}-r_n^k}{h_n},  \\
 \disp s &= & z- z_n,  \vphantom{\dfrac{1}{2}} \quad z\in [z_{n+1}, z_{n}].
 \ea
 \]
 
The piecewise--defined spline function $S(t_k,z)$ consists of the functions $S_n^k$ of the following on the interval $[z_n, z_{n+1}]$:
\be\label{speqn1}
\ba{lll}
\disp S_n^k(s)  & = &\disp  \dfrac{3h_ns^2-2s^3}{h_n^3}r_{n+1}^k + \dfrac{h_n^3-3h_ns^2+2s^3}{h_n^3}r_{n}^k \\
    && \qquad + \dfrac{s^2(s-h_n)}{h_n^2}d_{n+1}^k 
     + \dfrac{s(s-h_n)^2}{h_n^2}d_n^k,
 \ea
 \ee
for $0\le k \le K$.
Following the standard theory about cubic splines, we derive the following linear system
\be\label{speqn2}
\ba{l}
\disp \dfrac{3}{h_{n-1}}r_{n-1}^{k+1} + \Big( \dfrac{3}{h_n} -  \dfrac{3}{h_{n-1}} \Big)r_{n}^{k+1} - \dfrac{3}{h_n}r_{n+1}^{k+1}
 + h_n d_{n-1}^{\,k+1} \\
 \disp \qquad  + \,
2(h_{n-1}+h_n) d_{n}^{\,k+1} + h_{n-1} d_{n+1}^{\,k+1} = 0 \vphantom{\dfrac{1}{2}}, \qquad 
2\le n \le N\,.
\ea
\ee
Invoking the boundary conditions, the above system can be written as
\be\label{splin1}
  M_1 \,\vec{x} = \vec{b_1},
\ee
where $\vec{x}$ is given as in equation \eqref{unknown}, $\vec{b_1}$ is a column vector of length $N-1$ given by
\[
\disp \vec{b_1} = \begin{bmatrix}
 -h_2\alpha, 0, \cdots, 0, -h_{N}\beta\\
\end{bmatrix}^T\,,
\]
and $M_1=[M_{11},M_{12}]$ is an $(N-1)\times (2N)$ matrix with
\[
  \disp M_{11} = \begin{bmatrix}
  \dfrac{3}{h_{1}} & \Big(\dfrac{3}{h_2}-\dfrac{3}{h_{1}}\Big) & -\dfrac{3}{h_{2}} &  &   & \\
                   &  \dfrac{3}{h_{3}} & \Big(\dfrac{3}{h_2}-\dfrac{3}{h_{2 }}\Big) & -\dfrac{3}{h_{3}} &  & \\
                   &      \hphantom{\dfrac{3}{h_1}-\dfrac{3}{h_1}} \ddots &  \hphantom{\dfrac{3}{h_1}-\dfrac{3}{h_1}}\ddots &  \hphantom{\dfrac{3}{h_1}-\dfrac{3}{h_1}}\ddots \\               
                   &                 & \dfrac{3}{h_{N-1}} & \Big(\dfrac{3}{h_N}
                   -\dfrac{3}{h_{N-1}}\Big) & -\dfrac{3}{h_{N}} \\
  \end{bmatrix}_{(N-1)\times (N+1)\,,}
\]
and
\[
  \disp M_{12} = \begin{bmatrix}
    2(h_1+h_2) & h_1 &  & \\
       h_3     &  2(h_2+h_3) & h_2  & \\
      \hphantom{(h_2)}  \ddots      &\hphantom{(h_2)} \ddots & \hphantom{(h_2+h_3)}\ddots \\
               & h_{N-1}     & 2(h_{N-1}+h_N) & h_{N-2} \\
               &             &  h_N  & 2(h_{N-1}+h_N) 
  \end{bmatrix}_{(N-1)\times (N-1)\,.}
\]


On the other hand, from the equation \eqref{devolution}, we can derive the rest of the equations needed for solving $\vec{x}$. For convenience, we use
\[
    \disp (\mu_1)_n^k = \mu_1(r(t_k,z_n)), \quad (\mu_2)_n^k= \mu_2(r(t_k,z_n)),\quad
    Q_n^k = 1 + (r_z(t_k,z_n))^2,
\]
and
\[
    \disp \xi_n^k = \dfrac{1}{(\mu_1)_n^k}\dfrac{1}{Q_n^k} \tau,\quad
   \eta_n^k = - \Big(\dfrac{1}{(\mu_2)_n^k r_n^k}+ \overline{\Lambda}_n^{\,k}\sqrt{Q_n^k}\Big)\tau
\]
for $0\le k \le K$ and $1\le n\le N+1$. Then, the system \eqref{devolution} can be written as
\be\label{devolution2}
   r_n^{k+1} = (r_{zz})_n^{k+1}\xi_n^k  + \eta_n^k + r_n^{k}, \qquad 1\le n \le N+1.
\ee
Since from equation \eqref{speqn1}, we have
\[
  \ba{l}
   \disp r_n^{k+1} = \dfrac{6\delta_n^{k+1}- 2d_{n+1}^{\, k+1}- 4d_{n}^{\, k+1}}{h_n}, \qquad 1\le n\le N; \\
   \disp r_{n}^{k+1} = \dfrac{-6\delta_N^{k+1} + 4d_{N+1}^{\, k+1} + 2d_{N}^{\, k+1}}{h_N}, \qquad n=N+1,
\ea
\]
using the boundary conditions, we can write \eqref{devolution2}  as 
\be\label{splin2}
  M_2 \,\vec{x} = \vec{b_2},
\ee
where $\vec{b_2}$ is a column vector of length $N+1$ given by
\[
\disp \vec{b_2} = \begin{bmatrix}
 h_1(\eta_1^k+r_1^k)-4\alpha\xi_1^k\\
  h_2(\eta_2^k+r_2^k)\\
   \vdots\\
   h_{N-1}(\eta_{N-1}^k+r_{N-1}^k)\\
    h_{N}(\eta_{N}^k+r_{N}^k)- 2\beta\xi_N^k\\
  h_{N}(\eta_{N+1}^k+r_{N+1}^k)+ 4\beta\xi_{N+1}^k \\
\end{bmatrix}\,,
\]
and $M_2=[M_{21},M_{22}]$ is an $(N-1)\times (2N)$ matrix with
\[
  \disp M_{21} = \begin{bmatrix}
  h_1+\dfrac{6\xi_1^k}{h_1} & \hphantom{ \dfrac{6\xi_1^k}{h_1}} - \dfrac{6\xi_1^k}{h_1} &            \\
                            & \ddots                    & \ddots     \\
                             && h_N+\dfrac{6\xi_N^k}{h_N} & - \dfrac{6\xi_N^k}{h_N}  \\
                             & & - \dfrac{6\xi_{N+1}^k}{h_N} & h_N+\dfrac{6\xi_{N+1}^k}{h_N} \\
  \end{bmatrix}_{(N+1)\times (N+1)\,,}
\]
and
\[
  \disp M_{22} = \begin{bmatrix}
  2\xi_1^k &   &           &               &  \\
  4\xi_2^k & 2\xi_2^k      &               &  \\
           & \ddots        & \ddots        &  \\
           &               & 4\xi_{N-1}^k  & 2\xi_{N-1}^k  \\
           &               &               & 4\xi_{N}^k    \\
           &               &               & -2\xi_{N+1}^k   \\                           
  \end{bmatrix}_{(N+1)\times (N-1)\,.}
\]
Combining systems \eqref{splin1} and \eqref{splin2}, we finally obtain the linear system

\be\label{splin3}
  M\vec{x}= \vec{b}, 
\ee
where
\[
\disp M=
\begin{bmatrix}
  M_{11} & M_{12} \\
  M_{21} & M_{22}
\end{bmatrix}, \qquad \vec{b}= \begin{bmatrix}
 \vec{b_1}\\
 \vec{b_2}
\end{bmatrix}\, .
\]

As we can see, the system \eqref{splin3} is very sparse and can be solved by standard sparse solvers. We also would like to mention that the computation of $\overline{\Lambda}_n^{\,k}$ is done by using Gaussian quadratures.

\subsection{Numerical Experiments}

Some of the results from our numerical experiments will be presented below to illustrate the evolution the curvature flow that we have studied. In all of our simulations, we let $z\in [0,1]$.

\begin{example}\label{exp1}
In this experiment, $\epsilon = 0.2$, $T=3$, $N=500$, $\tau=10^{-4}$.
The initial profile of $r$ is chosen to be the cubic Hermite interpolant that satisfies
\[ \disp r(0,0)=0.7, \quad r(0,1)=0.4, \quad r_z(0,0)=r_z(0,1)=0. \]
\end{example}
 
The results are shown in Figure ~1:
\begin{itemize}
\item[(a)] $r(t,z)$, the profiles of the generating curve of the surface at different times are shown in plot (a). The initial profile is eventually evolves to a cylinder and remains thereafter. 
\item[(b)] The snapshots of the surfaces in three dimensional space at different times  are also shown in plot (b). 
\item[(c)] The history of the values of the energy functional $F$ (in equation \eqref{F}) is shown in plot (c), and it can be seen that the energy is decreasing as the surface evolves under equation \eqref{rt} until it remains nearly unchanged, which numerically implies that a minimum of the energy has been reached and the minimizer is corresponding to the the surface of a cylinder.
\item[(d)] The history of the values of the volume enclosed by the surface is shown in plot (d). Although it is seen that the volume is not preserved at the beginning of the evolution process, the error (relative to the initial volume) is within $0.007\%$. We think this is due to the error in the numerical approximations. First, the volumes are computed by using Gaussian quadrature. Secondly, the quantity $\overline{\Lambda}$ in equation \eqref{rt} is approximated by $\overline{\Lambda}_n^{\,k}$ as in our semi--implicit scheme, which numerically violates the law of volume preserving curvature flow (equation \eqref{rt}), unless in the later stage of the evolution, the surface nearly has constant anisotropic mean curvature. 
\end{itemize}

According to Theorem 5.1 of \cite{KP2006}, the threshold of stability for cylinders is 
\be\label{threshold}
\disp \dfrac{\mu_1(0)}{\mu_2(0)}\dfrac{1}{r^2}\le \dfrac{\pi^2}{h^2}\,.
\ee
For $\gamma = 1 + \epsilon\nu_3^2$ and $\epsilon = 0.2$, the cylinder is stable provided
\[ \disp r\ge \dfrac{\sqrt{1+2\epsilon}}{\pi}\approx 0.3766 .\]
This experiment numerically verifies the above stability analysis. In other similar experiments with larger volume fractions, we have also observed that no pinching has occurred.

\begin{figure}[h!]
\begin{center}
\begin{tabular}{cc}
   \includegraphics[scale=0.3]{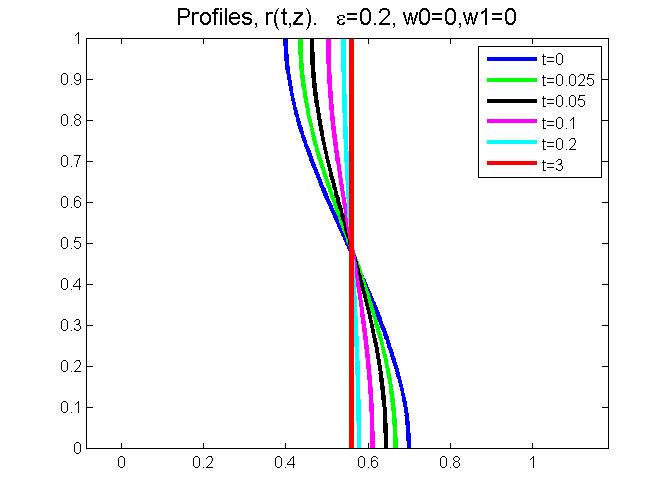} & \includegraphics[scale=0.3]{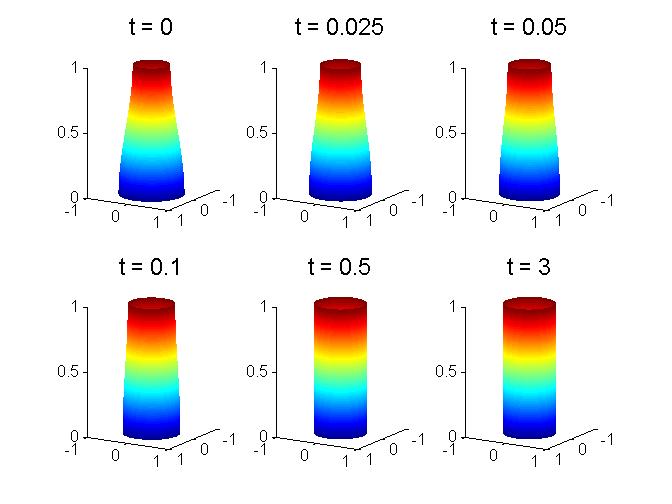}  \\
   (a) & (b)\\
   \\
   \includegraphics[scale=0.3]{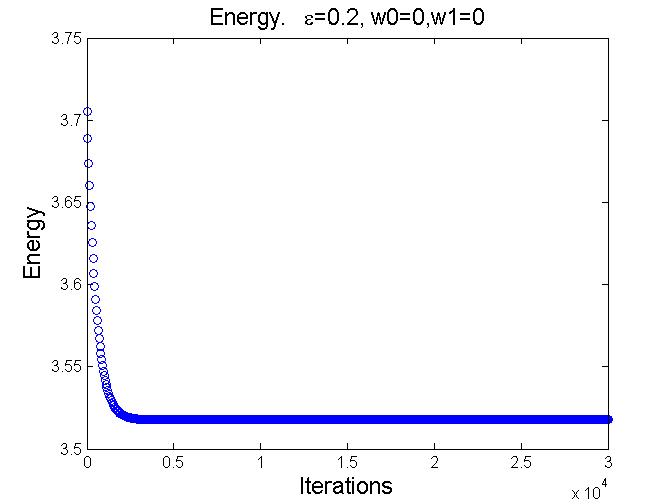} & \includegraphics[scale=0.3]{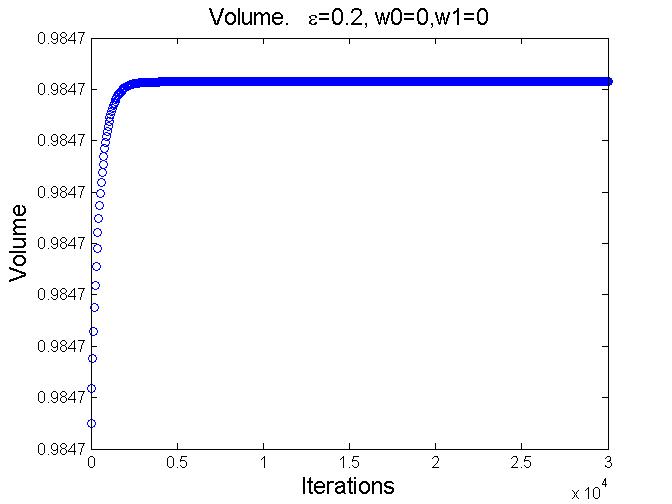} \\
   (c) & (d)\\
\end{tabular}
            
\end{center}
            \caption{Example \ref{exp1}. Snapshots of surface profiles, values of energy functional, and volumes. $\epsilon = 0.2$, $r(0,0)=0.7$, $r(0,1)=0.4$.}
            \label{fig1}
\end{figure}

\begin{example}\label{exp2}
 In this experiment, $\epsilon = 0.2$, $T=3$, $N=500$, $\tau=10^{-5}$.
The initial profile of $r$ is chosen to be the cubic Hermite interpolant that satisfies
\[ \disp r(0,0)=0.4, \quad r(0,1)=0.2, \quad r_z(0,0)=r_z(0,1)=0. \] 
\end{example}

\begin{figure}[h]
\begin{center}
\begin{tabular}{cc}
   \includegraphics[scale=0.3]{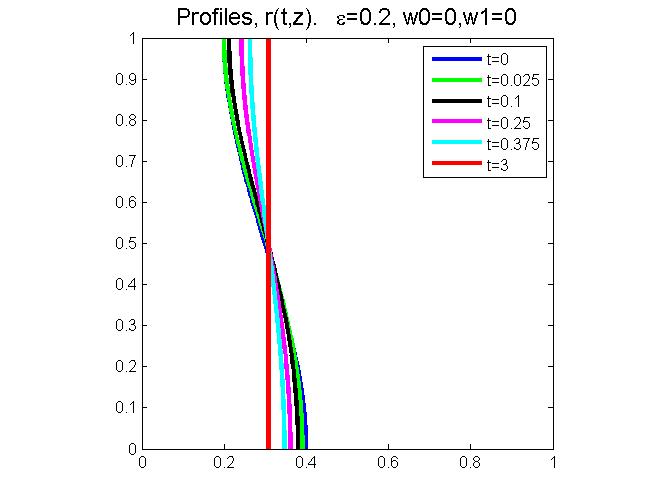} & \includegraphics[scale=0.3]{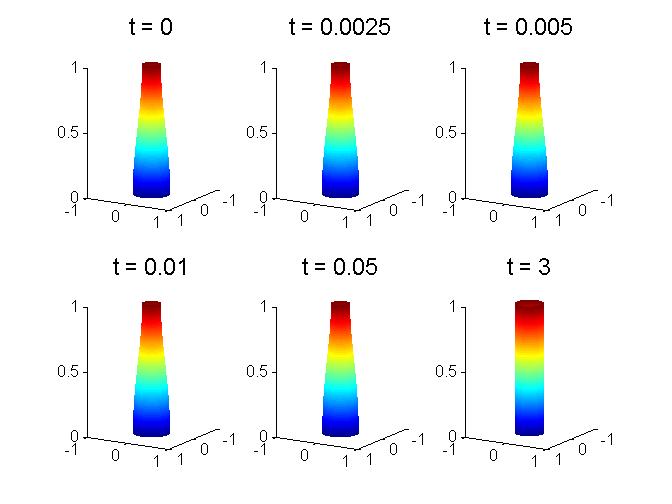}  \\
   (a) & (b)\\
   \\
   \includegraphics[scale=0.3]{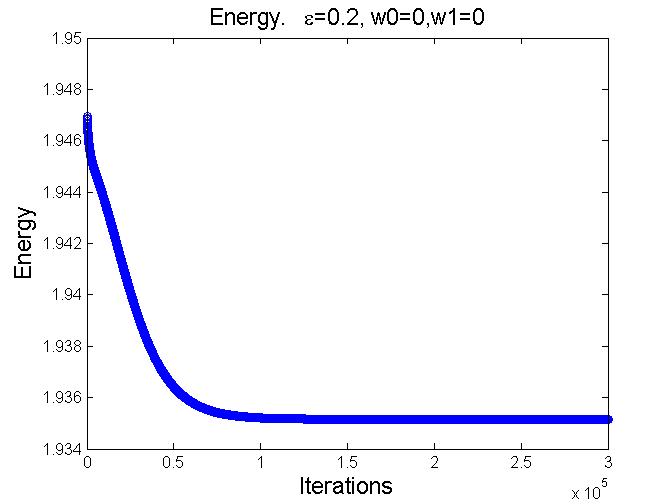} & \includegraphics[scale=0.3]{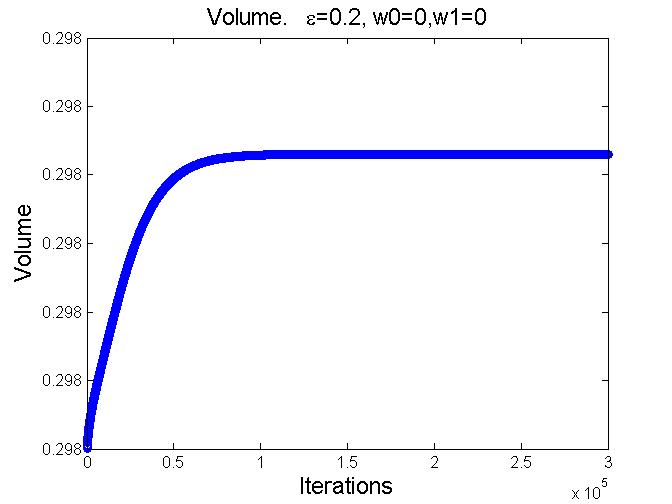} \\
   (c) & (d)
\end{tabular}        
\end{center}
     \caption{Example \ref{exp2}. Snapshots of surface profiles, values of energy functional, and volumes. $\epsilon = 0.2$, $r(0,0)=0.4$, $r(0,1)=0.2$.}
     \label{fig2}
\end{figure}

The results are shown in Figure ~\ref{fig2}. The initial volume is about $0.2980$. Comparing with $0.9847$, the volume in Example \ref{exp1}, this one is much smaller. However, the initial surface still evolves to a cylinder. Also, the errors in volumes (relative to the initial volume) are within $7\times 10^{-4}\,\%$.

\begin{figure}[h]
\begin{center}
\begin{tabular}{cc}
   \includegraphics[scale=0.3]{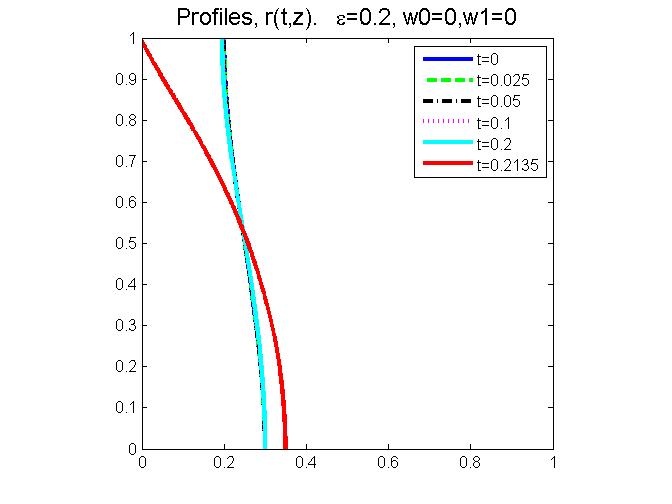} & \includegraphics[scale=0.3]{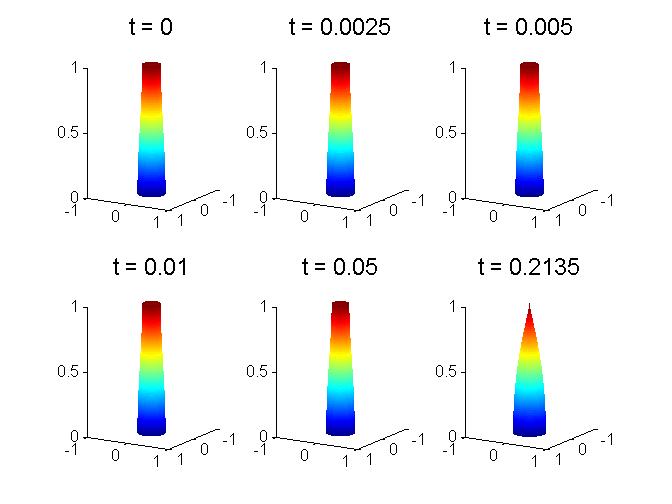}  \\
   (a) & (b)\\
   \\
   \includegraphics[scale=0.3]{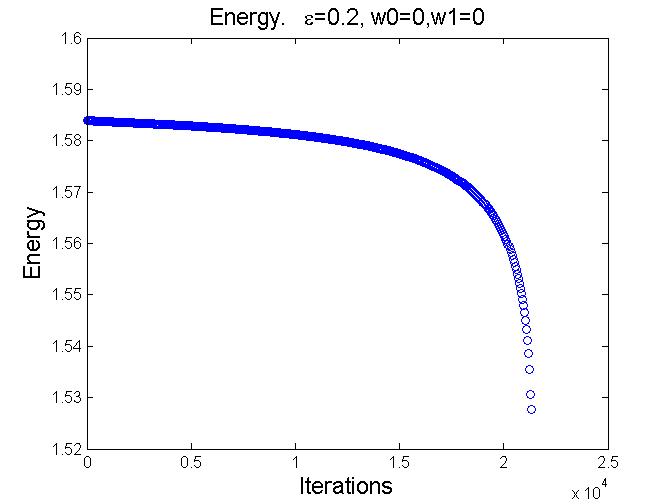} & \includegraphics[scale=0.3]{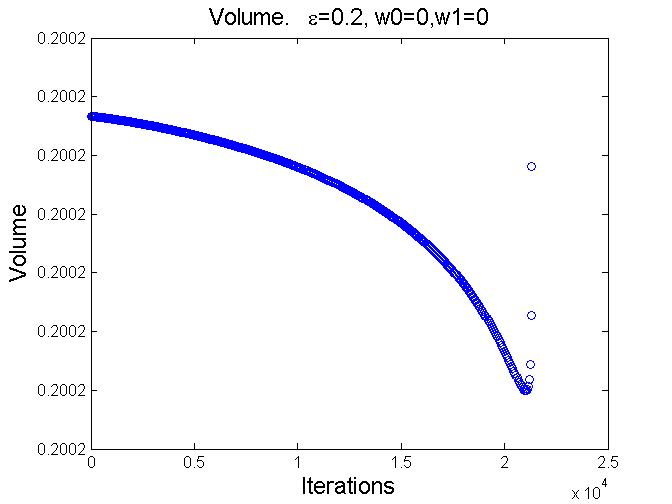} \\
   (c) & (d)
\end{tabular}        
\end{center}
     \caption{Example \ref{exp3}. Snapshots of surface profiles, values of energy functional, and volumes. $\epsilon = 0.2$, $r(0,0)=0.3$, $r(0,1)=0.2$.}
     \label{fig2a}
\end{figure}

\begin{example}\label{exp3}
In this experiment, $\epsilon = 0.2$, $T=3$, $N=500$, $\tau=10^{-5}$.
The initial profile of $r$ is chosen to be the cubic Hermite interpolant that satisfies
\[ \disp r(0,0)=0.3, \quad r(0,1)=0.2, \quad r_z(0,0)=r_z(0,1)=0. \] 
\end{example}
The only difference between this example and Example \ref{exp2} is the values of $r(0,0)$. The results are shown in Figure ~\ref{fig2a}. However, this time it can be seen that singularity occurs. 
 
\begin{figure}[h!]
\begin{center}
   \includegraphics[scale=0.3]{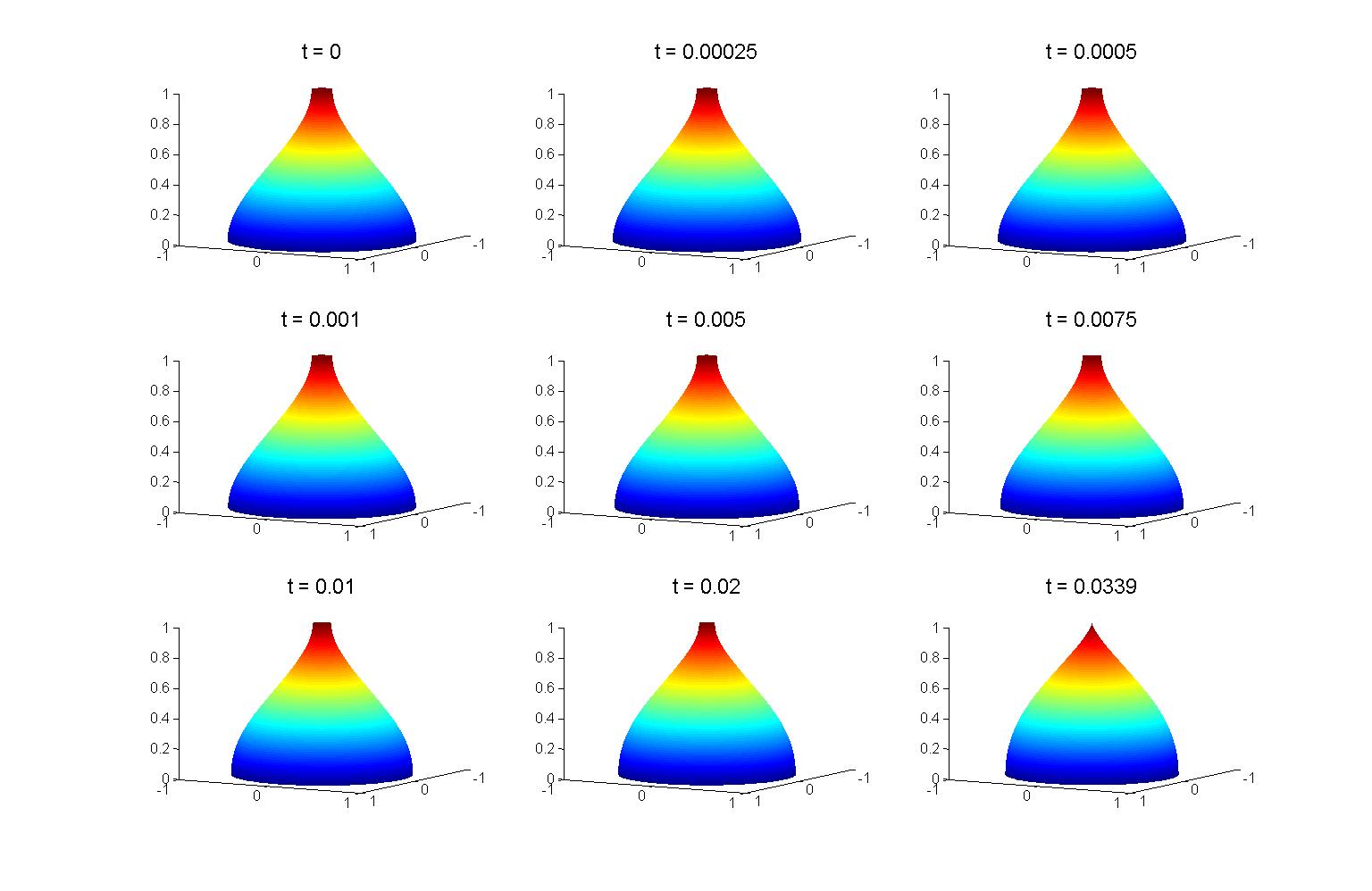}   
 \end{center}   
 \caption{Example \ref{exp4}. Snapshots of the 3D surfaces. $\epsilon = 0.2$, $r(0,0)=0.9$, $r(0,1)=0.1$.}
 \label{fig2b}
 \end{figure}

\begin{example}\label{exp4}
In this experiment, $\epsilon = 0.2$, $T=4$, $N=1000$, $\tau=10^{-6}$.
The initial profile of $r$ is chosen to be the cubic Hermite interpolant that satisfies
\[ \disp r(0,0)=0.9, \quad r(0,1)=0.1, \quad r_z(0,0)=r_z(0,1)=0. \]
The result is shown in Figure ~\ref{fig2b}. 
\end{example}

\begin{figure}[h!]
\begin{center}   
    \includegraphics[scale=0.3]{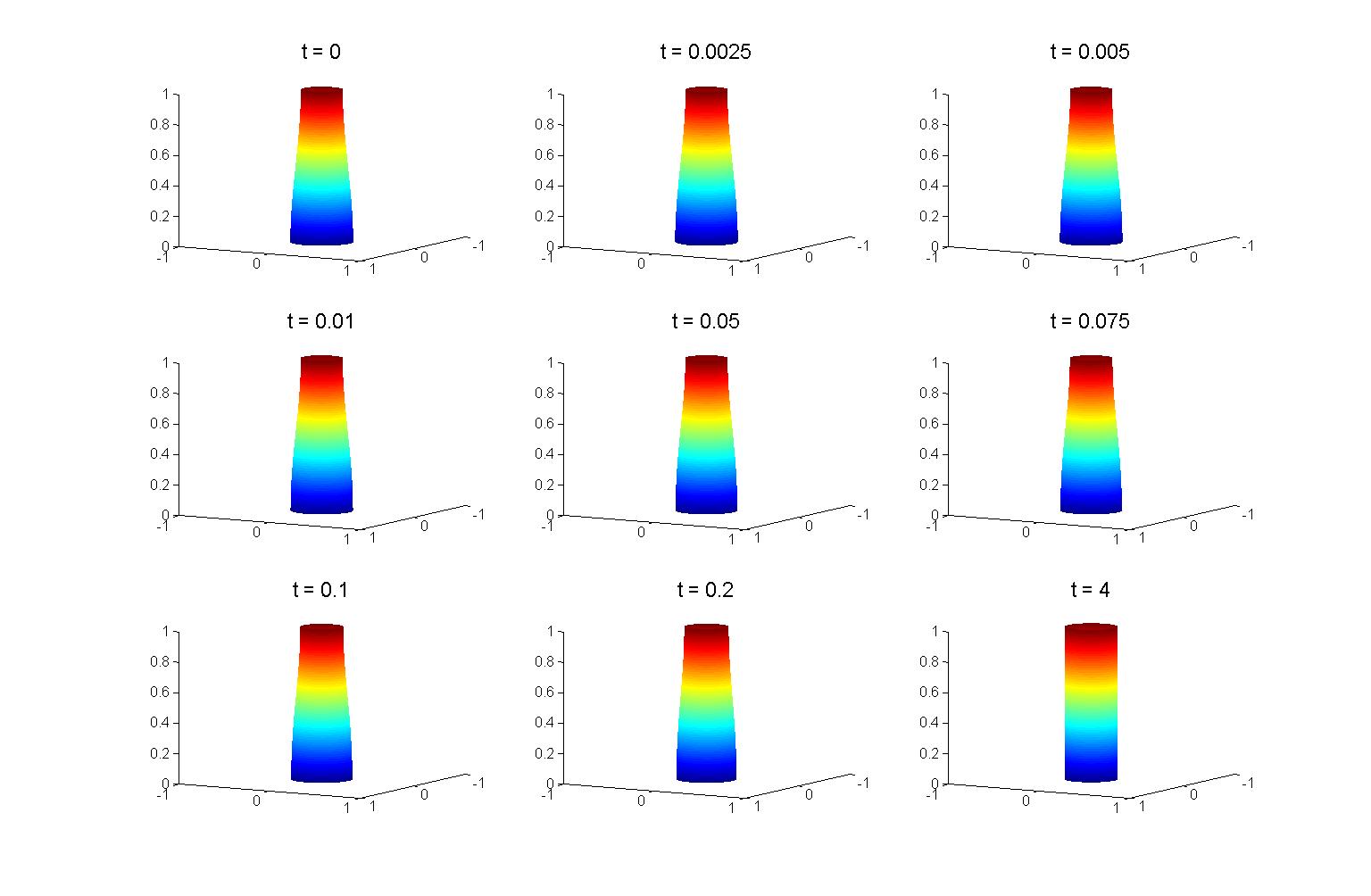}
 \end{center}   
 \caption{Example \ref{exp5}. Snapshots of the 3D surfaces. $\epsilon = 0.4$, $r(0,0)=0.3$, $r(0,1)=0.2$.}
 \label{fig2c}
 \end{figure}

\begin{example}\label{exp5}
In this experiment, $\epsilon = 0.4$, $T=4$, $N=500$, $\tau=10^{-5}$.
The initial profile of $r$ is chosen to be the cubic Hermite interpolant that satisfies
\[ \disp r(0,0)=0.3, \quad r(0,1)=0.2, \quad r_z(0,0)=r_z(0,1)=0. \] 
The results are shown in Figure ~\ref{fig2c}. 
\end{example}

Through these  examples \ref{exp2} to \ref{exp5}, we see that for small volume fractions, namely, when inequality \eqref{np} is not satisfied, the stability of cylinders is lost and singularity may develop, and therefore the flow may not exit for all time.

\begin{figure}[h]
\begin{center}
\begin{tabular}{cc}
   \includegraphics[scale=0.3]{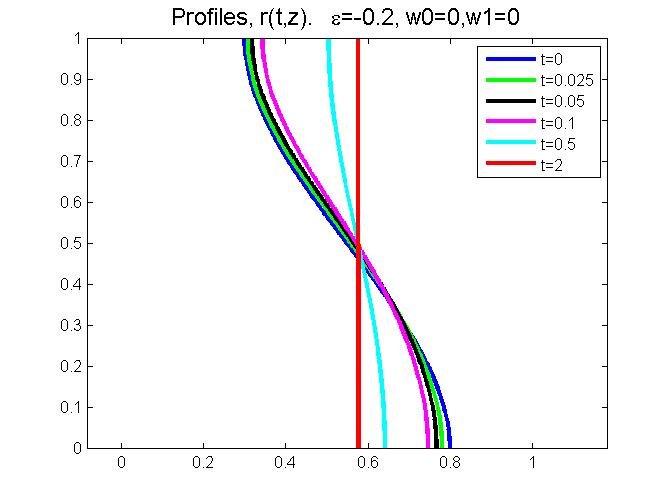} & \includegraphics[scale=0.3]{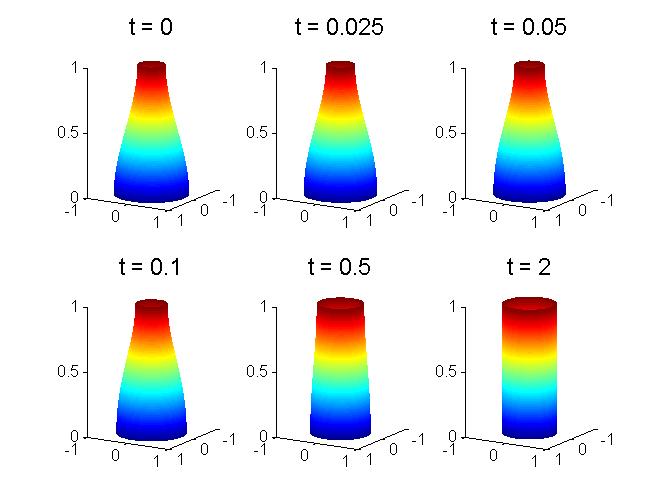}  \\
   (a) & (b) \\
   \\
   \includegraphics[scale=0.3]{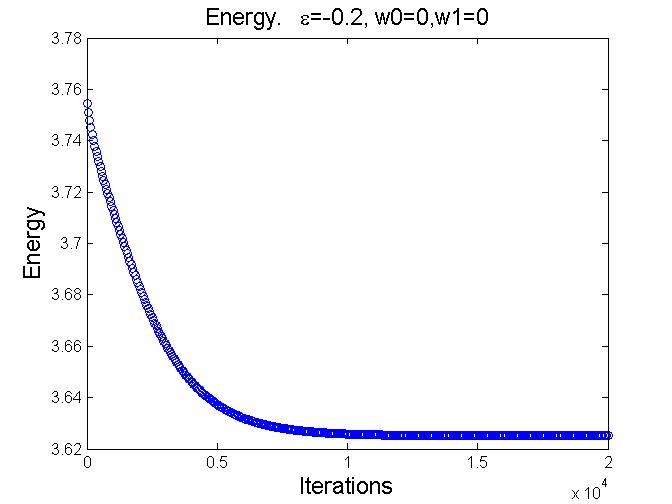} & \includegraphics[scale=0.3]{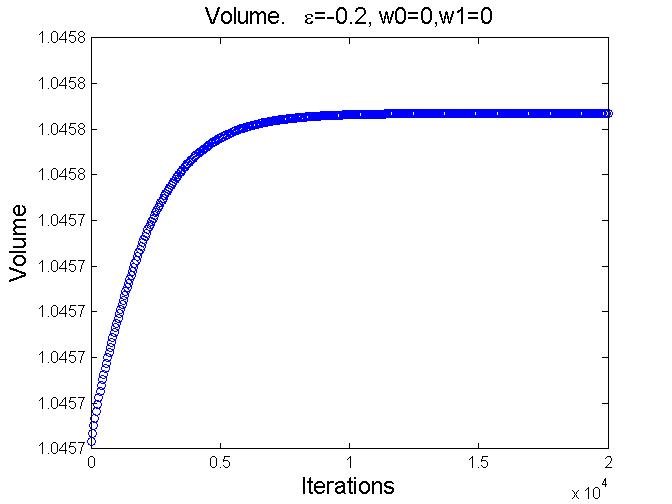} \\
   (c) & (d)
\end{tabular}
            
\end{center}
            \caption{Example \ref{exp6}. Snapshots of surface profiles, values of energy functional, and volumes. $\epsilon = -0.2$, $r(0,0)=0.8$, $r(0,1)=0.3$.}
            \label{fig3}
\end{figure}

\begin{figure}[h]
\begin{center}
\begin{tabular}{cc}
   \includegraphics[scale=0.25]{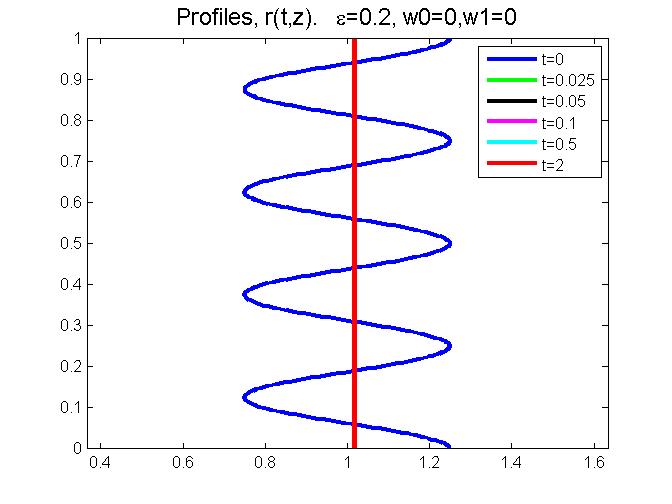} & \includegraphics[scale=0.25]{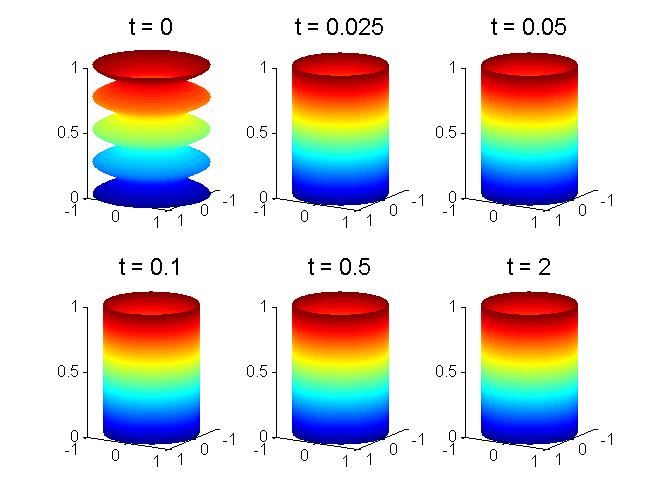}  \\
   (a) & (b) \\
   \\
   \includegraphics[scale=0.25]{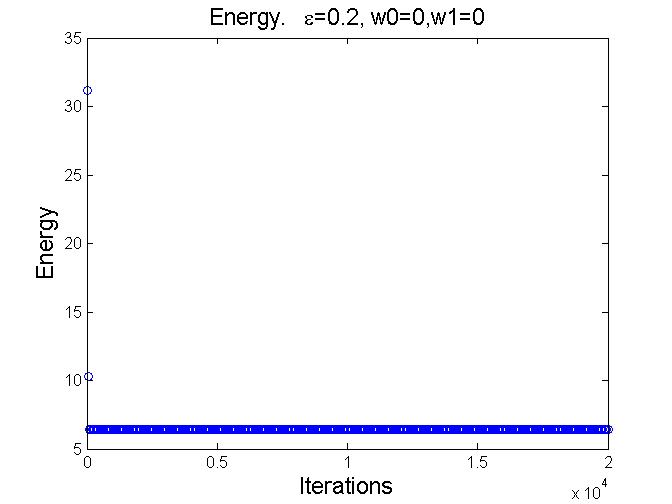} & \includegraphics[scale=0.25]{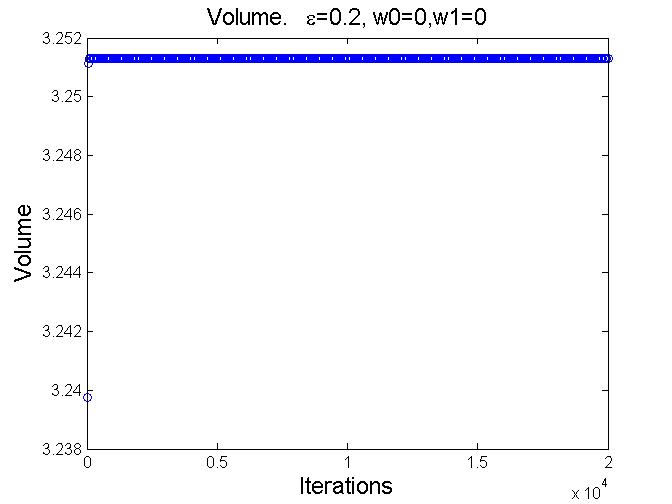} \\
   (c) & (d) 
\end{tabular}
            
\end{center}
      \caption{Snapshots of surface profiles, values of energy functional, and volumes. $\epsilon = 0.2$, $r(0,z)= 1 + \cos(8\pi z)/4$.}
      \label{fig4}
\end{figure}

We present two more examples: In Example \ref{exp6}, $\epsilon<0$, in Example \ref{exp7}, a different initial profile is used. The results are shown in Figure ~\ref{fig3} and Figure \ref{fig4}, respectively.

\begin{example}\label{exp6} In this experiment, $\epsilon = -0.2$, $T=2$, $N=500$, $\tau=10^{-4}$. The initial profile of $r$ is chosen to be the cubic Hermite interpolant that satisfies
\[ \disp r(0,0)=0.8, \quad r(0,1)=0.3, \quad r_z(0,0)=r_z(0,1)=0. \] 

\end{example}

\begin{example}\label{exp7} In this experiment, $\epsilon = 0.2$, the initial profile is given by
\[ \disp r(0,z)= 1 + \dfrac{1}{4} \cos(8\pi z). \]
\end{example}

\bigskip
\bigskip
\bigskip
\bigskip
\begin{flushleft}
Bennett P{\footnotesize ALMER}\\
Department of Mathematics\\
Idaho State University\\
Pocatello, ID 83209\\
U.S.A.\\
E-mail: palmbenn@isu.edu
\end{flushleft}
 \begin{flushleft}
Wenxiang Z{\footnotesize HU} \\
Department of Mathematics\\
Idaho State University\\
Pocatello, ID 83209\\
U.S.A.\\
E-mail: zhuwenx@isu.edu
\end{flushleft}

\end{document}